\documentclass{elsart}
\usepackage{amsfonts}

\usepackage{graphicx}
\usepackage{pstricks}
\usepackage{amsmath}
\usepackage{amsmath}
\usepackage{amsxtra}
\usepackage{graphicx}
\usepackage{pstricks}
\usepackage{amsmath}
\usepackage{amsmath}
\usepackage{amsxtra}

\newtheorem{theorem}{Theorem}[section]
\newtheorem{lemma}[theorem]{Lemma}

\newtheorem{example}[theorem]{Example}
\newtheorem{remark}[theorem]{Remark}

\numberwithin{equation}{section}
\begin{document}

\centerline{} 
\begin{frontmatter}


\title{Spectral pictures of $2$-variable weighted shifts\thanksref{label1}}
\thanks[label1]{Research partially supported by NSF Grants DMS-0099357 and
DMS-0400741.}


\author[authorlabel1]{Ra\'{u}l E. Curto},
\ead{rcurto@math.uiowa.edu}
\author[authorlabel2]{Jasang Yoon}
\ead{jyoon@iastate.edu}

\address[authorlabel1]{Department of Mathematics, The University of Iowa, Iowa City, Iowa
52242}
\address[authorlabel2]{Department of Mathematics, Iowa State University, Ames, Iowa 50011}


\medskip

\begin{abstract}
We study the spectral pictures of (jointly) hyponormal $2$-variable weighted
shifts with commuting subnormal components. \ By contrast with all known results
in the theory of
subnormal single and $2$-variable weighted shifts, we show that the Taylor
essential spectrum can be disconnected. \ We do this by obtaining a simple
sufficient condition that guarantees disconnectedness, based on the norms of
the horizontal slices of the shift. \ We also show that for every $k\geq 1$
there exists a $k$-hyponormal $2$-variable weighted shift whose 
horizontal and vertical slices have $1$- or $2$-atomic Berger measures, and 
whose Taylor
essential spectrum is disconnected.

\vskip 0.5\baselineskip

\noindent{\bf R\'esum\'e} \vskip 0.5\baselineskip \noindent
{\bf Les images spectrales de shifts pond\'{e}r\'{e}s \`a $2$-variables. }
Nous \'{e}tudions les images spectrales de shifts pond%
\'{e}r\'{e}s \`a deux variables et (conjointement) hyponormaux poss\'{e}dant
des composants commutants sousnormaux. \ \`{A} la diff\'{e}rence de tous les
r\'esultats connus dans la th\'{e}%
orie des shifts pond\'{e}r\'{e}s simples sousnormaux  \`a deux variables,
nous d\'{e}montrons que le spectre essentiel de Taylor peut \^{e}tre d\'{e}%
connect\'{e}. \ Nous faisons cela en obtenant une condition suffisante
simple qui garantit le caract\`{e}re d\'{e}connect\'{e} de ce spectre, bas%
\'{e}e sur les normes des sections horizontales du shift. \ Nous montrons 
\'{e}galement que pour chaque $k\geq 1$ il existe un shift pond\'{e}r\'{e}
$k$-hyponormal \`a deux variables dont le spectre essentiel  de Taylor
est d\'econnect\'{e}.
{\it Pour citer cet article~: A. Name1, A. Name2, C. R. Acad. Sci.
Paris, Ser. I 340 (2005).}

\end{abstract}
\end{frontmatter}


\section*{Version fran\c{c}aise abr\'eg\'ee}

Le Probl\`{e}me du Rel\`evement des Sousnormaux Commutants (PRSC) demande
quelles sont les conditions n\'{e}cessaires et suffisantes pour qu'une paire
commutante d'op\'{e}rateurs sousnormaux dans l'espace de Hilbert $\mathcal{H}$
admette une paire commutante d'extensions normales agissant sur un espace
Hilbert $\mathcal{K}\supseteq \mathcal{H}$. \ Dans nos recherches r\'{e}%
centes nous avons montr\'{e} que l'hyponormalit\'{e}, bien que n\'{e}%
cessaire, n'est pas une condition suffisante pour le rel\`evement. \ Nous
avons fait cela en consid\'{e}rant des shifts pond\'{e}r\'{e}s commutants en
deux variables, pour lesquels nous avons d\'{e}velopp\'{e} de nouvelles
techniques pour d\'{e}tecter leur hyponormalit\'{e} et sousnormalit\'{e}. \ 
\`{A} l'int\'{e}rieur de la classe $\mathfrak{H}_{0}$ de shifts pond\'{e}r%
\'{e}s de $2$-variables $\mathbf{T}\equiv (T_{1},T_{2})$ poss\'{e}dant des
composants sousnormaux $T_{1}$ et $T_{2}$, nous avons obtenu dans \cite%
{CuYo2} une nouvelle condition n\'{e}cessaire pour l'existence d'un
rel\`evement: les mesures de Berger des sections horizontales et verticales
doivent \^{e}tre ordonn\'{e}es de fa\c{c}on lin\'{e}aire par rapport \`{a}
la continuit\'{e} absolue. \ Plus r\'{e}cemment, dans \cite{CLY} nous avons
donn\'{e} une solution abstraite du PRSC, apr\`{e}s avoir d\'{e}montr\'{e}
une version multivari\'ee du crit\`{e}re de Bram-Halmos \cite{Con}.

Dans ce travail, nous d\'{e}veloppons des techniques et des outils nouveaux,
et nous les combinons avec ceux dans \cite{CuYo1}, \cite{CuYo2}, \cite{CuYo3}%
, \cite{CLY}, \cite{CuYa1} et \cite{CuYa2}, afin d'\'eclairer  la th\'{e}%
orie spectrale des shifts pond\'{e}r\'{e}s \`a $2$-variables. \ Comme il est
d\'eja bien connu, l'image spectrale d'un shift pond\'{e}r\'{e} unilat\'{e}%
ral et hyponormal $W_{\alpha }$ est facile \`{a} d\'{e}crire: le spectre est
un disque ferm\'{e} de rayon $\left\| W_{\alpha }\right\| $, le spectre
essentiel est un cercle de rayon $\left\| W_{\alpha }\right\| $, et l'indice
de Fredholm est $-1$ dans le disque ouvert. \ Ainsi, d'un point de vue 
spectral tout shift pond\'{e}r\'{e} unilat\'{e}ral hyponormal de norme $1$
est \'{e}quivalent au shift unilat\'{e}ral (non-pond\'{e}r\'{e}) $U_{+}$
(qui est aussi sousnormal). \ Pour les shifts pond\'{e}r\'{e}s \`a $2 $%
-variables la situation est tout \`{a} fait diff\'{e}rente, et une
description compl\`{e}te de l'image spectrale a \'{e}t\'{e} donn\'{e}e dans %
\cite{CuYa1} et \cite{CuYa2}, en utilisant la construction du groupo\"{\i}de
introduite dans \cite{Re} et \cite{MuR}, et raffin\'{e}e dans \cite{CuM}. \
La pr\'{e}sence de la mesure de Berger \'{e}tait essentielle dans l'\'{e}%
tude du comportement asymptotique des suites de poids, et a men\'{e} \`a des
r\'{e}sultats concrets sur les diverses parties du spectre de Taylor. \ Pour
les shifts pond\'{e}r\'{e}s \textit{hyponormaux} \`a $2$-variables $\mathbf{%
T} \in \mathfrak{H}_{0}$, pourtant, l'\'{e}tude des propri\'{e}t\'{e}s
spectrales exige des techniques ind\'{e}pendantes, puisque aucune mesure de
Berger n'est pr\'{e}sente. \ Dans ce qui suit, nous pr\'{e}sentons nombre de
r\'{e}sultats qui mettent en relief les diff\'{e}rences profondes entre le
cas en une variable et celui en deux variables. \ Dans le Th\'{e}or\`{e}me %
\ref{thmcompactper} nous exhibons, pour la premi\`{e}re fois, une condition
suffisante qui garantit le caract\`{e}re d\'{e}connect\'{e} du spectre
essentiel de Taylor, notamment, $\left\| W_{\alpha ^{(1)}}\right\| <\left\|
W_{\alpha ^{(0)}}\right\| $, o\`{u} $W_{\alpha ^{(j)}}$ d\'{e}note le $j$-%
\`{e}me niveau horizontal   de $T_{1}$; dans l'Exemple \ref{exampleBergman}
nous montrons que cette condition peut \^{e}tre pr\'{e}sente m\^{e}me dans
les paires hyponormales avec des sections $W_{\alpha ^{(j)}}\;(j\geq 1)$
mutuellement absolument continues. Nous am\'{e}liorons ce r\'{e}sultat en
montrant dans le Th\'{e}or\`eme \ref{important} qu'il est possible de former
davantage de composantes connexes du spectre essentiel de Taylor, tout en
pr\'{e}servant son caract\`{e}re hyponormal, si nous utilisons des shifts
pond\'{e}r\'{e}s de type Bergman dans chaque niveau horizontal de $T_{1}$. \
Ce fait est tout-\`a-fait \'{e}tonnant, \`{a} la lumi\`{e}re des r\'{e}%
sultats bien connus sur les shifts pond\'{e}r\'{e}s en une variable.
\ Au cours de notre analyse, nous d\'{e}montrons que pour les shifts pond%
\'{e}r\'{e}s sousnormaux en deux variables, les mesures de Berger des
sections horizontales $\{W_{\alpha ^{(j)}}\}_{j=1}^{\infty }$ sont toutes
mutuellement absolument continues et, par cons\'{e}quent, $\left\| W_{\alpha
^{(j)}}\right\| =\left\| W_{\alpha ^{(1)}}\right\| \;(j\geq 1)$. \ Cette
nouvelle condition n\'{e}cessaire est facilement calculable, et compl\`{e}te
la condition n\'{e}cessaire pr\'{e}c\'{e}dente qui se trouve dans %
\cite[Theorem 3.3]{CuYo2}. \ Ensuite, nous pr\'{e}sentons une application 
\`{a} la th\'{e}orie des op\'{e}rateurs \`{a} une variable (Exemple \ref%
{stair}). \ Finalement, le Th\'{e}or\`eme \ref{khypo of ex} montre que m\^{e}%
me si nous supposons que les mesures de Berger des shifts associ\'{e}s aux
lignes horizontales et colonnes verticales sont discr\`{e}tes, la $k$%
-hyponormalit\'{e} conjointe \cite{CLY} de $\mathbf{T}$ (pour un $k\geq 1$
arbitraire) n'est pas une condition suffisante pour forcer l'\'{e}galit\'{e} 
$\left| \sigma _{Te}(\mathbf{T})\right| =\partial \left| \sigma _{T}(\mathbf{%
T})\right| $. 

\section{Introduction}

The Lifting Problem for Commuting Subnormals (LPCS) asks for necessary and
sufficient conditions for a commuting pair of subnormal operators on Hilbert
space $\mathcal{H}$ to admit a commuting pair of normal extensions acting on
a Hilbert space $\mathcal{K}\supseteq \mathcal{H}$. \ In recent work we have
shown that (joint) hyponormality, while necessary, is not a sufficient
condition for lifting \cite{CuYo1}. \ We did this by appealing to commuting $%
2$-variable weighted shifts, for which we have developed new techniques to
detect their hyponormality and subnormality. \ Within the class $\mathfrak{H}%
_{0}$ of $2$-variable weighted shifts $\mathbf{T}\equiv (T_{1},T_{2})$ with
commuting subnormal components $T_{1}$ and $T_{2}$, we obtained in \cite%
{CuYo2} a new necessary condition for the existence of a lifting: the Berger
measures of horizontal and vertical slices must be linearly ordered with
respect to absolute continuity. \ More recently, we gave in \cite{CLY} an
abstract solution of LPCS, after proving a multivariable version of the
Bram-Halmos Criterion \cite{Con}.

In this paper we develop new tools and techniques, and combine them with
those in \cite{CuYo1}, \cite{CuYo2}, \cite{CuYo3}, \cite{CLY}, \cite{CuYa1}
and \cite{CuYa2}, to shed light on the spectral theory of $2$-variable
weighted shifts. \ As it is well known, the spectral picture of a hyponormal
unilateral weighted shift $W_{\alpha }$ is easy to describe: the spectrum is
a closed disk of radius $\left\| W_{\alpha }\right\| $, the essential
spectrum is the circle of radius $\left\| W_{\alpha }\right\| $, and the
Fredholm index is $-1$ in the open disk. \ Thus, from a spectral perspective
all norm-one hyponormal unilateral weighted shifts are equivalent to the
unilateral (unweighted) shift $U_{+}$ (which is also subnormal). \ For
subnormal $2$-variable weighted shifts the situation is quite different; a
complete description of the spectral picture, in the case when the
intersection of the boundary of the Taylor spectrum and each coordinate
plane is a circle, was given in \cite{CuYa1} and \cite{CuYa2}. This was done
using the groupoid machinery introduced in \cite{Re} and \cite{MuR}, and
refined in \cite{CuM}. \ The presence of the Berger measure was essential in
the study of the asymptotic behavior of the weight sequences, and led to
concrete results about the various parts of the Taylor spectrum. \ For 
\textit{hyponormal} $2$-variable weighted shifts $\mathbf{T} \in \mathfrak{H}%
_{0}$, however, the study of the spectral properties requires independent
techniques, since no Berger measure is present.

In what follows, we present a number of results which highlight the deep
differences between the single variable case and the $2$-variable case. \ In
Theorem \ref{thmcompactper} we exhibit, for the first time, a sufficient
condition that guarantees the \textit{disconnectedness} of the Taylor
essential spectrum, namely, $\left\| W_{\alpha ^{(1)}}\right\| <\left\|
W_{\alpha ^{(0)}}\right\| $, where $W_{\alpha ^{(j)}}$ denotes the $j$-th
horizontal slice of $T_{1}$; in Example \ref{exampleBergman} we show that
this condition can even be present in hyponormal pairs with mutually
absolutely continuous horizontal slices $W_{\alpha ^{(j)}}\;(j\geq 1)$. \ We
improve this by showing in Theorem \ref{important} that more connected
components of the Taylor essential spectrum can be formed, still preserving
hyponormality, if we use Bergman-like weighted shifts on each horizontal
level of $T_{1}$. \ This fact is quite surprising, in view of the well known
one-variable results. \ Along the way we prove that, for subnormal $2$%
-variable weighted shifts, the Berger measures of the horizontal slices $%
\{W_{\alpha ^{(j)}}\}_{j=1}^{\infty }$ are all mutually absolutely
continuous and, as a result, $\left\| W_{\alpha ^{(j)}}\right\| =\left\|
W_{\alpha ^{(1)}}\right\| \;($all $j\geq 1)$. \ This new necessary condition
is easily computable, and complements the previous necessary condition found
in \cite[Theorem 3.3]{CuYo2}. \ Next, we present an application to single
variable operator theory. \ It is well known that for $T$ a hyponormal
operator on $\mathcal{H}$, $r((T-\lambda )^{-1})=\frac{1}{dist(\lambda
,\sigma (T))}$ ($\lambda \notin \sigma (T)$), where $r$ denotes spectral
radius and $dist$ denotes distance. \ If we substitute the left spectrum for
the spectrum, the result is far from obvious; we actually show in Example %
\ref{stair} that for $T$ hyponormal, the equality $\left\| [(T-\lambda
)^{(\ell )}]^{-1}\right\| =\frac{1}{dist(\lambda ,\sigma _{\ell }(T))}$ may
fail, where $\sigma _{\ell }$ denotes the left spectrum and $[(T-\lambda
)^{(\ell )}]^{-1}:=[(T-\lambda )^{\ast }(T-\lambda )]^{-1/2}$ is the
canonical left inverse. \ Finally, Theorem \ref{khypo of ex} shows that even
if we assume that the Berger measures of the shifts associated to horizontal
rows and vertical columns are discrete, joint $k$-hyponormality \cite{CLY}
of $\mathbf{T}$ (for an arbitrary $k\geq 1$) is not sufficient to force the
equality $\left| \sigma _{Te}(\mathbf{T})\right| =\partial \left| \sigma
_{T}(\mathbf{T})\right| $. \ (For compact sets $K\subseteq \mathbb{C}^{2}$
and $L\subseteq \mathbb{R}^{2}$, $\left| K\right| :=\{(\left| z_{1}\right|
,\left| z_{2}\right| ):(z_{1},z_{2})\in K\}$ and $\partial L$ denotes the
outer boundary of $L$, that is, the boundary of the connected component of $%
\mathbb{R}^{2}\backslash L$.)

\section{Main Results}

We begin by listing three basic results which are needed in the proofs of
Theorems \ref{thmcompactper}, \ref{important} and \ref{khypo of ex}, and
Example \ref{exampleBergman}.

\begin{lemma}
\label{lem1}$(i)$ (\cite{Cu3}, \cite{Cu1}) \ Let $\mathcal{H}_{1}$ and $\mathcal{H}_{2}$
be Hilbert spaces, and let $A_{i}\in \mathcal{L}(\mathcal{H}_{1}),$ $%
C_{i}\in \mathcal{L}(\mathcal{H}_{2})$ and $B_{i}\in \mathcal{L}(\mathcal{H}%
_{1},\mathcal{H}_{2}),(i=1,\cdots ,n)$ be such that $\left( 
\begin{array}{cc}
\mathbf{A} & \mathbf{0} \\ 
\mathbf{B} & \mathbf{C}%
\end{array}%
\right) :=\left( \left( 
\begin{array}{cc}
A_{1} & 0 \\ 
B_{1} & C_{1}%
\end{array}%
\right) ,\ldots ,\left( 
\begin{array}{cc}
A_{n} & 0 \\ 
B_{n} & C_{n}%
\end{array}%
\right) \right) $ is commuting. Assume that $\mathbf{A}$ and $\left( 
\begin{array}{cc}
\mathbf{A} & \mathbf{0} \\ 
\mathbf{B} & \mathbf{C}%
\end{array}%
\right) $ are Taylor invertible. Then $\mathbf{C}$ is Taylor invertible.%
\newline
$(ii)$ (\cite{CuFi}, \cite{Fia} and \cite{Cu3}) \ For $\mathbf{A}$ and $\mathbf{B}$ two
commuting $n$-tuples of bounded operators on Hilbert space, we have $\sigma
_{T}(\mathbf{A}\otimes I,I\otimes \mathbf{B})=\sigma _{T}(\mathbf{A})\times \sigma _{T}(\mathbf{B})$, $\sigma
_{\ell }(\mathbf{A}\otimes I,I\otimes \mathbf{B})=\sigma _{\ell }(\mathbf{A})\times \sigma _{\ell }(\mathbf{B})$
and $\sigma _{r}(\mathbf{A}\otimes I,I\otimes \mathbf{B})=\sigma _{r}(\mathbf{A})\times \sigma _{r}(\mathbf{B}),$
where $\sigma _{\ell }$ and $\sigma _{r}$ denote the left and right spectra,
respectively.\newline
$(iii)$ (\cite{CuYo2}) \ Let $\mu $ be the Berger measure of a subnormal $2$%
-variable weighted shift, and for $j\geq 0$ let $\xi _{j}$ be the Berger
measure of the associated $j$-th horizontal $1$-variable weighted shift $%
W_{\alpha ^{(j)}}$. Then\ $\xi _{j}=\mu _{j}^{X}$ (the marginal measure of $%
\mu _{j}$), where $d\mu _{j}(s,t):=\frac{1}{\gamma _{0j}}t^{j}d\mu (s,t)$;
more precisely, $d\xi _{j}(s)=\{\frac{1}{\gamma _{0j}}\int_{Y}t^{j}\;d\Phi
_{s}(t)\}\;d\mu ^{X}(s),$ where $d\mu (s,t)\equiv d\Phi _{s}(t)\;d\mu ^{X}(s)
$ is the disintegration of $\mu $ by vertical slices. A similar result holds
for the Berger measure $\eta _{i}$ of the associated $i$-th vertical $1$%
-variable weighted shifts $W_{\beta ^{(i)}}\;(i\geq 0)$.
\end{lemma}

\begin{theorem}
\label{thmcompactper}Let $\mathbf{T}$ be the $2$-variable weighted shift
given by Figure \ref{compactperturbation}, let $W_{\alpha
^{(j)}}:=shift(\alpha _{0j},\alpha _{1j},\cdots )$ $(j\geq 0)
$, $W_{\beta }:=shift(\beta _{0},\beta _{1},\cdots )$, and let $B:=\beta
_{0}diag(1,\frac{\alpha _{01}}{\alpha _{00}},\cdots )$. \ Assume that $%
\left\| W_{\alpha ^{(j)}}\right\| =\left\| W_{\alpha ^{(1)}}\right\|
<\left\| W_{\alpha ^{(0)}}\right\| \;($all $j\geq 2$). Then $B$ is a
compact operator, and $\sigma _{T}(\mathbf{T})=(\left\| W_{\alpha
^{(1)}}\right\| \cdot \overline{\mathbb{D}}\times \left\| W_{\beta }\right\|
\cdot \overline{\mathbb{D}})\cup (\left\| W_{\alpha ^{(0)}}\right\| \cdot 
\overline{\mathbb{D}}\times \mathbb{\{}0\})$, $\sigma _{Te}(\mathbf{T}%
)=(\left\| W_{\alpha ^{(1)}}\right\| \cdot \overline{\mathbb{D}}\times
\left\| W_{\beta }\right\| \cdot \mathbb{T})\cup (\left\| W_{\alpha
^{(1)}}\right\| \cdot \mathbb{T}\times \left\| W_{\beta }\right\| \cdot 
\overline{\mathbb{D}})\cup (\left\| W_{\alpha ^{(0)}}\right\| \cdot \mathbb{T%
}\times \mathbb{\{}0\}).$ \ In particular, $\sigma _{Te}(\mathbf{T})\neq
\partial \sigma _{T}(\mathbf{T})$. $\ $(Here $\overline{\mathbb{D}}$ denotes
the closure of the open unit disk $\mathbb{D}$, $\mathbb{T}$ the unit
circle, and $\partial K:=\{(z_{1},z_{2})\in \mathbb{C}^{2}:(\left|
z_{1}\right| ,\left| z_{2}\right| )\in \partial \left| K\right| \}$).
\end{theorem}

\begin{theorem}
\label{necessary cond}For $i,j\geq 1$, let $W_{\alpha ^{(j)}}$ $%
( $resp. $W_{\beta ^{(i)}}:=shift(\beta _{i0},\beta _{i1},\cdots
)) $ be the $j$-th horizontal slice (resp. $i$-th vertical slice) of
a subnormal $2$-variable weighted shift. Then $\left\Vert W_{\alpha
^{(j)}}\right\Vert =\left\Vert W_{\alpha ^{(1)}}\right\Vert $ and $%
\left\Vert W_{\beta ^{(i)}}\right\Vert =\left\Vert W_{\beta
^{(1)}}\right\Vert $.
\end{theorem}

\noindent \textbf{Proof.} \ As a consequence of Lemma \ref{lem1}(iii), we
see that $\xi _{j}\approx \xi _{1}$ and $\eta _{i}\approx \eta
_{1} $ for all $i,j\geq 1$. \ It follows that $supp\xi _{j}=supp\xi _{1}$
and $supp\eta _{i}=supp\eta _{1}$ for all $i,j\geq 1$. Since the norm of a
subnormal unilateral weighted shift always equals the supremum of the
support of its Berger measure \cite{Con}, the result follows. \qed

\setlength{\unitlength}{1mm} \psset{unit=1mm} 
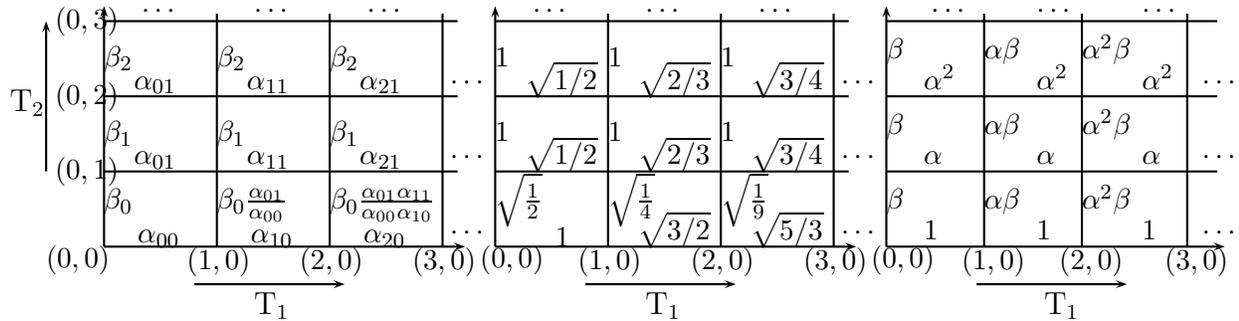
\begin{figure}[th]
\begin{center}
\begin{picture}(162,36.5)

\psline{->}(10,6)(58,6)
\psline(10,16)(57,16)
\psline(10,26)(57,26)
\psline(10,36)(57,36)
\psline{->}(10,6)(10,38)
\psline(25,6)(25,37)
\psline(40,6)(40,37)
\psline(55,6)(55,37)

\put(2.3,3.2){\footnotesize{$(0,0)$}}
\put(21,3){\footnotesize{$(1,0)$}}
\put(36,3){\footnotesize{$(2,0)$}}
\put(51,3){\footnotesize{$(3,0)$}}

\put(14.5,6.5){\footnotesize{$\alpha_{00}$}}
\put(29.5,6.5){\footnotesize{$\alpha_{10}$}}
\put(44.5,6.5){\footnotesize{$\alpha_{20}$}}
\put(56,7){\footnotesize{$\cdots$}}

\put(14,17){\footnotesize{$\alpha_{01}$}}
\put(29,17){\footnotesize{$\alpha_{11}$}}
\put(44,17){\footnotesize{$\alpha_{21}$}}
\put(56,17){\footnotesize{$\cdots$}}

\put(14,27){\footnotesize{$\alpha_{01}$}}
\put(29,27){\footnotesize{$\alpha_{11}$}}
\put(44,27){\footnotesize{$\alpha_{21}$}}
\put(56,27){\footnotesize{$\cdots$}}

\put(15,36.5){\footnotesize{$\cdots$}}
\put(30,36.5){\footnotesize{$\cdots$}}
\put(45,36.5){\footnotesize{$\cdots$}}

\psline{->}(22,1)(42,1)
\put(30,-3){$\rm{T}_1$}
\psline{->}(2.3, 16)(2.3,36)
\put(-2.5,24){$\rm{T}_2$}

\put(3.5,15){\footnotesize{$(0,1)$}}
\put(3.5,25){\footnotesize{$(0,2)$}}
\put(3.5,35){\footnotesize{$(0,3)$}}

\put(10,11){\footnotesize{$\beta_{0}$}}
\put(10,20){\footnotesize{$\beta_{1}$}}
\put(10,30){\footnotesize{$\beta_{2}$}}

\put(25,11){\footnotesize{$\beta_{0}\frac{\alpha_{01}}{\alpha_{00}}$}}
\put(25,20){\footnotesize{$\beta_{1}$}}
\put(25,30){\footnotesize{$\beta_{2}$}}

\put(40,11){\footnotesize{$\beta_{0}\frac{\alpha_{01}\alpha_{11}}{\alpha_{00}\alpha_{10}}$}}
\put(40,20){\footnotesize{$\beta_{1}$}}
\put(40,30){\footnotesize{$\beta_{2}$}}


\psline{->}(62,6)(110,6)
\psline(62,16)(109,16)
\psline(62,26)(109,26)
\psline(62,36)(109,36)
\psline{->}(62,6)(62,38)
\psline(77,6)(77,37)
\psline(92,6)(92,37)
\psline(107,6)(107,37)

\put(60,3.2){\footnotesize{$(0,0)$}}
\put(73,3){\footnotesize{$(1,0)$}}
\put(88,3){\footnotesize{$(2,0)$}}
\put(103,3){\footnotesize{$(3,0)$}}

\put(70,6){\footnotesize{$1$}}
\put(81,7){\footnotesize{$\sqrt{3/2}$}}
\put(96,7){\footnotesize{$\sqrt{5/3}$}}
\put(108,7){\footnotesize{$\cdots$}}

\put(66,17){\footnotesize{$\sqrt{1/2}$}}
\put(81,17){\footnotesize{$\sqrt{2/3}$}}
\put(96,17){\footnotesize{$\sqrt{3/4}$}}
\put(108,17){\footnotesize{$\cdots$}}

\put(66,27){\footnotesize{$\sqrt{1/2}$}}
\put(81,27){\footnotesize{$\sqrt{2/3}$}}
\put(96,27){\footnotesize{$\sqrt{3/4}$}}
\put(108,27){\footnotesize{$\cdots$}}

\put(67,36.5){\footnotesize{$\cdots$}}
\put(82,36.5){\footnotesize{$\cdots$}}
\put(97,36.5){\footnotesize{$\cdots$}}

\psline{->}(74,1)(94,1)
\put(82,-3){$\rm{T}_1$}

\put(62,11){\footnotesize{$\sqrt{\frac{1}{2}}$}}
\put(62,20){\footnotesize{$1$}}
\put(62,30){\footnotesize{$1$}}

\put(77,11){\footnotesize{$\sqrt{\frac{1}{4}}$}}
\put(77,20){\footnotesize{$1$}}
\put(77,30){\footnotesize{$1$}}

\put(92,11){\footnotesize{$\sqrt{\frac{1}{9}}$}}
\put(92,20){\footnotesize{$1$}}
\put(92,30){\footnotesize{$1$}}


\psline{->}(114,6)(159,6)
\psline(114,16)(158,16)
\psline(114,26)(158,26)
\psline(114,36)(158,36)
\psline{->}(114,6)(114,38)
\psline(127,6)(127,37)
\psline(140,6)(140,37)
\psline(154,6)(154,37)

\put(112,3.2){\footnotesize{$(0,0)$}}
\put(124,3){\footnotesize{$(1,0)$}}
\put(137,3){\footnotesize{$(2,0)$}}
\put(151,3){\footnotesize{$(3,0)$}}

\put(119,7){\footnotesize{$1$}}
\put(134,7){\footnotesize{$1$}}
\put(148,7){\footnotesize{$1$}}
\put(156,7){\footnotesize{$\cdots$}}

\put(119,17){\footnotesize{$\alpha$}}
\put(134,17){\footnotesize{$\alpha$}}
\put(148,17){\footnotesize{$\alpha$}}
\put(156,17){\footnotesize{$\cdots$}}

\put(119,27){\footnotesize{$\alpha^{2}$}}
\put(134,27){\footnotesize{$\alpha^{2}$}}
\put(148,27){\footnotesize{$\alpha^{2}$}}
\put(156,27){\footnotesize{$\cdots$}}

\put(120,36.5){\footnotesize{$\cdots$}}
\put(134,36.5){\footnotesize{$\cdots$}}
\put(148,36.5){\footnotesize{$\cdots$}}

\psline{->}(126,1)(146,1)
\put(135,-3){$\rm{T}_1$}

\put(114,11){\footnotesize{$\beta$}}
\put(114,21){\footnotesize{$\beta$}}
\put(114,31){\footnotesize{$\beta$}}

\put(127,11){\footnotesize{$\alpha\beta$}}
\put(127,21){\footnotesize{$\alpha\beta$}}
\put(127,31){\footnotesize{$\alpha\beta$}}

\put(140,11){\footnotesize{$\alpha^{2}\beta$}}
\put(140,21){\footnotesize{$\alpha^{2}\beta$}}
\put(140,31){\footnotesize{$\alpha^{2}\beta$}}

\end{picture}
\end{center}
\par
\medskip
\caption{Weight diagrams of the 2-variable weighted shifts in Theorem \ref%
{thmcompactper} and Examples \ref{exampleBergman} and \ref{exof1atom},
respectively}
\label{compactperturbation}
\end{figure}

Example \ref{exampleBergman} below provides a concrete instance of Theorem %
\ref{thmcompactper}, with all horizontal rows admitting continuous Berger
measures. We first need some notation and a few definitions. Let $\mu $ be a
Reinhardt measure on $\mathbb{C}^{n}.$ The set of bounded point evaluations
for $\mu $ is $b.p.e(\mu )=\{\mathbf{\lambda }\in \mathbb{C}%
^{n}:p\rightarrow p(\mathbf{\lambda }),$ $p\in \mathbb{C[}z]$, extends to a
bounded point evaluation from $P^{2}(\mu )$ to $\mathbb{C}.\}$. The kernel
function associated with $\mu $ is $\emph{k}(z,w)\equiv \emph{k}^{(\mu
)}(z,w):=\sum_{\alpha \in \mathbb{Z}_{+}^{n}}\frac{z^{\alpha }\overline{w}%
^{\alpha }}{\left\| z^{\alpha }\right\| ^{2}}$, and the set of convergence
of $\emph{k}$ is $\mathcal{C}(\emph{k}):=\{\mathbf{\lambda }\in \mathbb{C}%
^{n}:\emph{k}(\mathbf{\lambda },\mathbf{\lambda })<\infty \}$.

\begin{example}
\label{exampleBergman}Let $\mathbf{T}\equiv (T_{1},T_{2})$ be the $2$
-variable weighted shift whose weight diagram is given in Figure \ref%
{compactperturbation}. Then $\mathbf{T} \in \mathfrak{H}_{0}$, $\mathbf{T}$
is hyponormal, $\xi _{j}\approx \xi _{1}\;(j\geq 1)$, $\eta
_{i}\approx \eta _{1}\;(i\geq 0)$ and $supp\eta
_{i}=\{0,1\}\;(i\geq 0)$. \ However, $\mathbf{T}$ is not subnormal.
Moreover $\sigma _{T}(\mathbf{T)=}\sigma _{r}(\mathbf{T})=\overline{\mathcal{%
C}(\emph{k})}=(\overline{\mathbb{D}}\times \overline{\mathbb{D}})\cup (\sqrt{%
2}\cdot \overline{\mathbb{D}}\times \{0\})$, $\sigma _{Te}(\mathbf{T}%
)=\sigma _{re}(\mathbf{T})=(\sqrt{2}\mathbb{T},0)\cup ((\mathbb{T}\times 
\overline{\mathbb{D}})\cup (\overline{\mathbb{D}}\times \mathbb{T}))$, and $%
\partial \sigma _{T}(\mathbf{T})\neq \sigma _{Te}(\mathbf{T})$ (cf. Figure %
\ref{compact}).
\end{example}

Theorem \ref{important} below shows that, for every $k\geq 1$, it is
possible to create $k$ connected components in the $\sigma _{Te}$ of a $2$%
-variable weighted shift $\mathbf{T}$, while maintaining the hyponormality
of $\mathbf{T} \in \mathfrak{H}_{0}$. We recall that for $\ell \geq 1$, the
Bergman-like weighted shift on $\ell ^{2}(\mathbb{Z}_{+})$ is $B_{+}^{(\ell
)}:=shift(\{\sqrt{\ell -\frac{1}{n+2}}:n\geq 0\})$; in particular, $%
B_{+}^{(1)}\equiv B_{+}$ is the Bergman shift. Note that $\left\|
B_{+}^{(\ell )}\right\| =\sqrt{\ell }\geq 1$, $B_{+}^{(\ell )}$ is
subnormal, with Berger measure $\xi _{\ell }$ $(\ell \geq 1),$ and $d\xi
_{1}(s)=ds$ on $[0,1]$ and $d\xi _{2}(s)=$$\frac{sds}{\pi \sqrt{2s-s^{2}}}$
on $[0,2]$ (cf. \cite{CPY}, \cite{CuYo3})$.$ \ The following result is a
variation of \cite[Theorem 3.14]{CuYo3}.

\begin{theorem}
\label{important2}For every $k\geq 1$ there exists $(i)$ a family of $%
\{B_{+}^{(\ell _{j})}\}_{j=0}^{k-1}$ of Bergman-like weighted shifts, and\ $%
(ii)$ a subnormal weighted shift $W_{\beta }\equiv shift(\beta _{0},\beta
_{1},\cdots )$ $($with $\beta _{n}<\beta _{n+1}$ all $n\geq 0)$, such that
the commuting $2$-variable weighted shift $\mathbf{T}$ with a weight diagram
whose first $k$ rows are $B_{+}^{(\ell _{0})},\cdots ,B_{+}^{(\ell _{k-1})}$%
, whose remaining rows are all equal to $U_{+}$, and whose $0$-th column is
given by $W_{\beta }$, is hyponormal.
\end{theorem}

\begin{theorem}
\label{important}Let $\mathbf{T}$ be a hyponormal $2$-variable weighted
shift satisfying the hypotheses in Theorem \ref{important2}, and let $%
b^{(\ell _{j})}:=\left\| B_{+}^{(\ell _{j})}\right\| \;(j=0,\cdots ,k-1)$
and $c:=\left\| W_{\beta }\right\| $. Then $\sigma _{T}(\mathbf{T})=(%
\overline{\mathbb{D}}\times c\overline{\mathbb{D}})\cup (b^{(\ell _{0})}%
\overline{\mathbb{D}}\times \{0\})$, $\sigma _{Te}(\mathbf{T})=[\overline{%
\mathbb{D}}\times \{c\}]\cup \lbrack \{1\}\times c\overline{\mathbb{D}}]\cup
\lbrack (b^{(\ell _{0})}\mathbb{T}\cup \cdots \cup b^{(\ell _{k-1})}\mathbb{%
T)\times \{}0\}]$ (cf. Figure \ref{compact}).
\end{theorem}

Example \ref{exof1atom} shows that, if we don't insist that $\mathbf{T} \in  
\mathfrak{H}_{0}$ be hyponormal, the Taylor essential spectrum can consist
of infinitely many circles converging to a single point. We see in
particular that two commuting subnormals do not necessarily give rise to a
hyponormal pair.

\begin{example}
\label{exof1atom}Let $\mathbf{T}$ be the $2$-variable weighted shift whose
weight diagram is given by Figure \ref{compactperturbation}, where $\alpha
<\beta \leq 1$. \ Then $\mathbf{T} \in \mathfrak{H}_{0}$, with $1$-atomic
Berger measures for all horizontal and vertical slices. \ However, $\mathbf{T%
}$ is not hyponormal. \ Moreover, $\sigma _{T}(\mathbf{T})$ has empty
interior; in fact, $\sigma _{T}(\mathbf{T})=\sigma _{r}(\mathbf{T})=%
[(\left| z_{1}\right| \leq 1)\times
\{0\}]\bigcup \ [\{0\}\times (\left| z_{2}\right| \leq \beta )],$ $\sigma
_{\ell }(\mathbf{T})=\sigma _{\ell e}(\mathbf{T})=\sigma _{re}(\mathbf{T}),$ 
$\sigma _{Te}(\mathbf{T})=\{(0,0)\}\bigcup \ \{[\bigcup_{k=0}^{\infty }(\left|
z_{1}\right| =\alpha ^{k})]\times \{0\}\}\bigcup \ \{\{0\}\times \lbrack
\bigcup \ _{\ell =0}^{\infty }(\left| z_{2}\right| =\beta \alpha ^{\ell })]\}$
and $\sigma _{\ell }(\mathbf{T})=\sigma _{Te}(\mathbf{T})$ (cf. Figure \ref%
{compact}).
\end{example}

We now turn to a simple application of $2$-variable weighted shifts to a
single variable problem. Example \ref{stair} shows that there exists a
hyponormal $T\in \mathcal{L}(\mathcal{H})$ with $\left\| (T-\lambda )^{(\ell
)^{-1}}\right\| \neq \frac{1}{dist(\lambda ,\sigma _{\ell }(T))}$ for some $%
\lambda \notin \sigma _{\ell }(T)$.

\begin{example}
\label{stair}Let $\mathbf{T}$ be the $2$-variable weighted shift whose
weight diagram is given by Figure \ref{compact}, where $a<1$. \ Then (i) $%
\mathbf{T}$\ is not hyponormal, and $\sigma _{\ell }(\mathbf{T})=\sigma
_{\ell e}(\mathbf{T})=(a\mathbb{T\times T})\cup (\mathbb{T\times }a\mathbb{T)%
}$; (ii) $T_{1}$ is hyponormal, $\sigma _{\ell }(T_{1})=\sigma _{\ell
e}(T_{1})=a\mathbb{T}\cup \mathbb{T}$ and at least one horizontal slice $%
W_{\alpha ^{(j)}}$ of $T_{1}$ does not satisfy the identity $\left\|
(W_{\alpha ^{(j)}}-\lambda )^{(\ell )^{-1}}\right\| =\frac{1}{dist(\lambda
,\sigma _{\ell }(W_{\alpha ^{(j)}}))}$ for all $\lambda \in \sigma _{\ell
}(W_{\alpha ^{(j)}})$.

\setlength{\unitlength}{1mm} \psset{unit=1mm} 
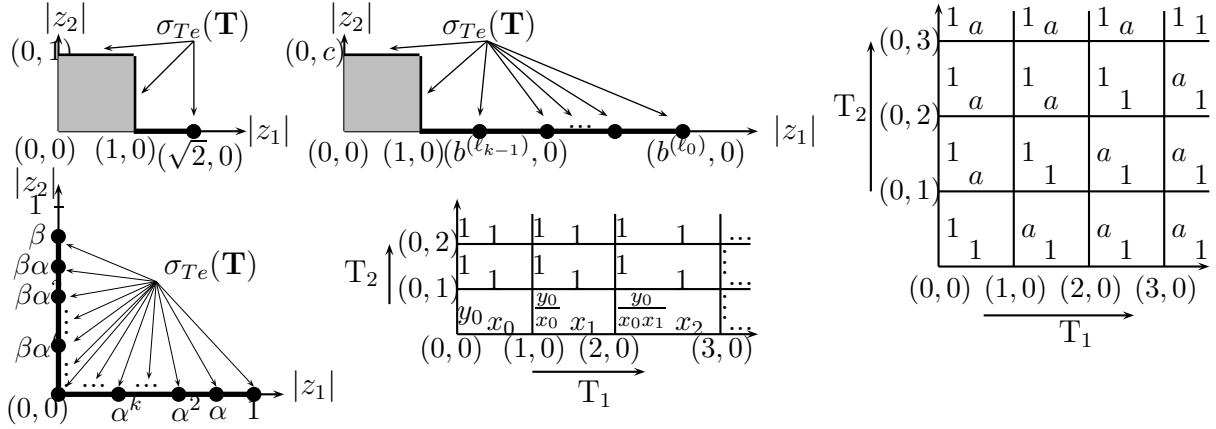
\begin{figure}[th]
\begin{center}
\begin{picture}(162,49.5)


\psline{->}(7,35)(7,48)
\psline{->}(7,35)(32,35)

\put(1,31.5){\footnotesize{$(0,0)$}}
\put(11.7,31.5){\footnotesize{$(1,0)$}}
\put(0.5,44.5){\footnotesize{$(0,1)$}}

\put(20,30.5){\footnotesize{$(\sqrt{2},0)$}}
\put(25,35){\pscircle*(0,0){1}}

\put(5,49){$|z_2|$}
\put(32,34){$|z_1|$}

\psline[linewidth=2pt](7,45)(17,45)
\psline[linewidth=2pt](17,45)(17,35)
\psline[linewidth=2pt](17,35)(25,35)

\pspolygon*[linecolor=lightgray](7,35)(7,45)(17,45)(17,35)
\psline[linewidth=0.5pt]{<-}(25,37)(25,47)
\psline[linewidth=0.5pt]{<-}(18,40)(25,47)
\psline[linewidth=0.5pt]{<-}(13,46)(25,47)
\put(20,48){$\sigma _{Te}(\mathbf{T})$}


\psline{->}(45,35)(45,48)
\psline{->}(45,35)(100,35)

\put(40,31){\footnotesize{$(0,0)$}}
\put(50,31){\footnotesize{$(1,0)$}}
\put(37,44.5){\footnotesize{$(0,c)$}}

\put(58,31){\footnotesize{$(b^{({\ell _{k-1}})},0)$}}
\put(63,35){\pscircle*(0,0){1}}
\put(72,35){\pscircle*(0,0){1}}
\put(74.7,35.5){$...$}
\put(81,35){\pscircle*(0,0){1}}
\put(90,35){\pscircle*(0,0){1}}
\put(85,31){\footnotesize{$(b^{({\ell _{0}})},0)$}}

\put(43,49){$|z_2|$}
\put(102,34){$|z_1|$}

\psline[linewidth=2pt](45,45)(55,45)
\psline[linewidth=2pt](55,45)(55,35)
\psline[linewidth=2pt](55,35)(90,35)

\pspolygon*[linecolor=lightgray](45,35)(45,45)(55,45)(55,35)

\psline[linewidth=0.5pt]{<-}(51,46)(64,47)
\psline[linewidth=0.5pt]{<-}(56,40)(64,47)
\psline[linewidth=0.5pt]{<-}(63,37)(64,47)
\psline[linewidth=0.5pt]{<-}(71.5,37)(64,47)
\psline[linewidth=0.5pt]{<-}(75.3,37)(64,47)
\psline[linewidth=0.5pt]{<-}(80,37)(64,47)
\psline[linewidth=0.5pt]{<-}(88,37)(64,47)
\put(58,48){$\sigma _{Te}(\mathbf{T})$}


\psline{->}(124,17)(160,17)
\psline(124,27)(160,27)
\psline(124,37)(160,37)
\psline(124,47)(160,47)
\psline{->}(124,17)(124,52)
\psline(134,17)(134,51)
\psline(144,17)(144,51)
\psline(154,17)(154,51)

\put(116,26){\footnotesize{$(0,1)$}}
\put(116,36){\footnotesize{$(0,2)$}}
\put(116,46){\footnotesize{$(0,3)$}}
\put(120,13){\footnotesize{$(0,0)$}}
\put(130,13){\footnotesize{$(1,0)$}}
\put(140,13){\footnotesize{$(2,0)$}}
\put(150,13){\footnotesize{$(3,0)$}}

\put(128,18){\footnotesize{$1$}}
\put(138,18){\footnotesize{$1$}}
\put(148,18){\footnotesize{$1$}}
\put(158,18){\footnotesize{$1$}}

\put(128,28){\footnotesize{$a$}}
\put(138,28){\footnotesize{$1$}}
\put(148,28){\footnotesize{$1$}}
\put(158,28){\footnotesize{$1$}}

\put(128,38){\footnotesize{$a$}}
\put(138,38){\footnotesize{$a$}}
\put(148,38){\footnotesize{$1$}}
\put(158,38){\footnotesize{$1$}}

\put(128,48){\footnotesize{$a$}}
\put(138,48){\footnotesize{$a$}}
\put(148,48){\footnotesize{$a$}}
\put(158,48){\footnotesize{$1$}}

\psline{->}(130,10.5)(150,10.5)
\put(140,7){$\rm{T}_1$}
\psline{->}(115, 27)(115,47)
\put(110,37){$\rm{T}_2$}

\put(125,21){\footnotesize{$1$}}
\put(125,31){\footnotesize{$1$}}
\put(125,41){\footnotesize{$1$}}
\put(125,49){\footnotesize{$1$}}

\put(135,21){\footnotesize{$a$}}
\put(135,31){\footnotesize{$1$}}
\put(135,41){\footnotesize{$1$}}
\put(135,49){\footnotesize{$1$}}

\put(145,21){\footnotesize{$a$}}
\put(145,31){\footnotesize{$a$}}
\put(145,41){\footnotesize{$1$}}
\put(145,49){\footnotesize{$1$}}

\put(155,21){\footnotesize{$a$}}
\put(155,31){\footnotesize{$a$}}
\put(155,41){\footnotesize{$a$}}
\put(155,49){\footnotesize{$1$}}


\psline{->}(7,0)(7,28)
\psline{->}(7,0)(37,0)
\psline[linewidth=2pt](7,0)(7,21)
\psline[linewidth=2pt](7,0)(33,0)

\put(0,-3){\footnotesize{$(0,0)$}}

\put(32,-3.3){\footnotesize{$1$}}
\put(27,-3.3){\footnotesize{$\alpha$}}
\put(22,-3.5){\footnotesize{$\alpha^{2}$}}
\put(17,1){$...$}
\put(14,-3.5){\footnotesize{$\alpha^{k}$}}
\put(10,1){$...$}

\put(33,0){\pscircle*(0,0){1}}
\put(28,0){\pscircle*(0,0){1}}
\put(23,0){\pscircle*(0,0){1}}
\put(15,0){\pscircle*(0,0){1}}
\put(7,0){\pscircle*(0,0){1}}

\put(3,24){\footnotesize{$1$}}

\put(7,21){\pscircle*(0,0){1}}
\put(3,20){\footnotesize{$\beta$}}

\put(7,17){\pscircle*(0,0){1}}
\put(1,16){\footnotesize{$\beta\alpha$}}

\put(7,13){\pscircle*(0,0){1}}
\put(1.2,12){\footnotesize{$\beta\alpha^{2}$}}
\put(7.5,8){$\vdots $}
\put(7,6.5){\pscircle*(0,0){1}}
\put(7.5,1){$\vdots $}
\put(1.2,5){\footnotesize{$\beta\alpha^{\ell}$}}

\put(1,27){$|z_2|$}
\put(38,0){$|z_1|$}
\put(6.2,24.8){$\_$}

\put(21,16){$\sigma _{Te}(\mathbf{T})$}
\psline[linewidth=0.3pt]{<-}(33,1)(20,15)
\psline[linewidth=0.3pt]{<-}(28,1)(20,15)
\psline[linewidth=0.3pt]{<-}(23,1)(20,15)
\psline[linewidth=0.3pt]{<-}(19,2)(20,15)
\psline[linewidth=0.3pt]{<-}(15,1)(20,15)
\psline[linewidth=0.3pt]{<-}(12,2)(20,15)
\psline[linewidth=0.3pt]{<-}(8,1)(20,15)
\psline[linewidth=0.3pt]{<-}(8,20)(20,15)
\psline[linewidth=0.3pt]{<-}(8,16.5)(20,15)
\psline[linewidth=0.3pt]{<-}(8.5,13)(20,15)
\psline[linewidth=0.3pt]{<-}(9,10)(20,15)
\psline[linewidth=0.3pt]{<-}(9,7)(20,15)
\psline[linewidth=0.3pt]{<-}(9,4)(20,15)


\psline{->}(60,8)(100,8)
\psline(60,14)(99,14)
\psline(60,20)(99,20)
\psline{->}(60,8)(60,26)
\psline(70,8)(70,24)
\psline(81,8)(81,24)
\psline(95,8)(95,24)

\put(52,13){\footnotesize{$(0,1)$}}
\put(52,19){\footnotesize{$(0,2)$}}

\put(55,5){\footnotesize{$(0,0)$}}
\put(66,5){\footnotesize{$(1,0)$}}
\put(76,5){\footnotesize{$(2,0)$}}
\put(91,5){\footnotesize{$(3,0)$}}

\put(64,8.5){\footnotesize{$x_{0}$}}
\put(75,8.5){\footnotesize{$x_{1}$}}
\put(89,8.5){\footnotesize{$x_{2}$}}
\put(96,8.5){$...$}

\put(64,14){\footnotesize{$1$}}
\put(75,14){\footnotesize{$1$}}
\put(89,14){\footnotesize{$1$}}
\put(96,14.5){$...$}

\put(64,20){\footnotesize{$1$}}
\put(75,20){\footnotesize{$1$}}
\put(89,20){\footnotesize{$1$}}
\put(96,20.5){$...$}

\psline{->}(70,3)(85,3)
\put(76,-1){$\rm{T}_1$}
\psline{->}(51,12)(51,20)
\put(45,15){$\rm{T}_2$}

\put(60,10){\footnotesize{$y_{0}$}}
\put(60,16){\footnotesize{$1$}}
\put(60,21){\footnotesize{$1$}}

\put(70,10.5){\footnotesize{$\frac{y_{0}}{x_{0}}$}}
\put(70,16){\footnotesize{$1$}}
\put(70,21){\footnotesize{$1$}}

\put(81,10.5){\footnotesize{$\frac{y_{0}}{x_{0}x_{1}}$}}
\put(81,16){\footnotesize{$1$}}
\put(81,21){\footnotesize{$1$}}

\put(95,10){$\vdots$}
\put(95,16){$\vdots$}

\end{picture}
\end{center}
\par
\medskip
\caption{(top left and center) Spectral pictures in Example \ref%
{exampleBergman} and Theorem \ref{important}, respectively; (bottom left)
spectral picture in Example \ref{exof1atom}; (top right and bottom center)
weight diagrams of the $2$-variable weighted shifts in Example \ref{stair}
and Theorem \ref{khypo of ex}, respectively}
\label{compact}
\end{figure}
\end{example}

The next result, Theorem \ref{khypo of ex}, shows that $k$-hyponormality ($%
k\geq 1$) is not sufficient to guarantee $\partial \sigma _{T}(\mathbf{T}%
)=\sigma _{Te}(\mathbf{T})$. We first recall that if $\mathbf{T}$ is a $2$%
-variable weighted shift with weight sequences $\alpha \equiv \{\alpha _{%
\mathbf{m}}\}$ and $\beta \equiv \{\beta _{\mathbf{m}}\}$, then $\mathbf{T}$
is $k$-hyponormal $\Leftrightarrow M_{\mathbf{m}}(k)\geq 0$ for all $\mathbf{%
m}\in \mathbb{Z}_{+}^{2}$, where $M_{\mathbf{m}}(k)$ is the $(k+1)\times
(k+1)$ matrix of moments beginning at $\mathbf{m}$ (\cite[Theorem 3.1]{CLY}).

\begin{lemma}
\label{mini lemma}$(i)$ For $\kappa >1$, let $x\equiv
\{x_{n}\}_{n=0}^{\infty }$ where $x_{0}:=\sqrt{\frac{1+\kappa }{2}}$ and $%
x_{n}:=\sqrt{\frac{1+\kappa ^{n+1}}{1+\kappa ^{n}}}\;(n\geq 1\}.$ Then $%
W_{\kappa }:=shift(x_{0},x_{1},\cdots )$ is subnormal, with Berger measure $%
\xi _{\kappa }=\frac{1}{2}(\delta _{1}+\delta _{\kappa })$.\newline
$(ii)$ Given $k\geq 1$ and $m_{1}\geq 0$, consider $M(m_{1}):=\left( 
\begin{array}{ccc}
\gamma _{m_{1}}(W_{\kappa }) & \cdots & \gamma _{m_{1}+k}(W_{\kappa }) \\ 
\vdots & \ddots & \vdots \\ 
\gamma _{m_{1}+k}(W_{\kappa }) & \cdots & \gamma _{m_{1}+2k}(W_{\kappa })%
\end{array}%
\right) $ and $\mathbf{1}:=\left( 
\begin{array}{ccc}
1 & \cdots & 1 \\ 
\vdots & \ddots & \vdots \\ 
1 & \cdots & 1%
\end{array}%
\right) $. If $0<y_{0}^{2}\leq \frac{1}{2}$ then $M(m_{1})-y_{0}^{2}\cdot
\mathbf{1}\geq 0$, where $\mathbf{1}$ is a $(k+1)\times (k+1)$ matrix.
\end{lemma}

\begin{theorem}
\label{khypo of ex}Let $\mathbf{T}$ be the $2$-variable weighted shift whose
weight diagram is given in Figure \ref{compact}. \ Given $k\geq 1$, $\kappa
>1$ and $W_{\kappa }$, there exists $0<y_{0}\leq \sqrt{\frac{1}{2}}$ which
makes $\mathbf{T}$ $k$-hyponormal. \ Moreover, $\sigma _{T}(\mathbf{T)=}%
\sigma _{r}(\mathbf{T})=(\overline{\mathbb{D%
}}\times \overline{\mathbb{D}})\cup (\sqrt{\kappa }\cdot \overline{\mathbb{D}%
}\times \{0\})$, $\sigma _{Te}(\mathbf{T})=\sigma _{re}(\mathbf{T})=(\sqrt{%
\kappa }\mathbb{T} \times 0)\cup ((\mathbb{T}\times \overline{\mathbb{D}})\cup (%
\overline{\mathbb{D}}\times \mathbb{T}))$, and $\partial \sigma _{T}(\mathbf{%
T})\neq \sigma _{Te}(\mathbf{T})$.
\end{theorem}

\noindent \textbf{Sketch of Proof.} \ Observe that the restriction of $%
\mathbf{T}$ to $\mathbf{\vee }\{e_{(m_{1},m_{2})}:m_{2}\geq 1\}$ is
unitarily equivalent to $(I\otimes U_{+},U_{+}\otimes I)$. \ We then apply (%
\cite[Theorem 3.1]{CLY}) to $\mathbf{m:}=(m_{1},0)$, with $m_{1}\geq 0$. \
We have $M_{(m_{1},0)}(k)\geq 0$ $($all $k\geq 1)\Leftrightarrow
M(m_{1})-y_{0}^{2}\cdot 1\geq 0$. \ We then use the moments associated with $%
\mathbf{T}$, matrix row operations, Choleski's algorithm, and a choice of $%
y_{0}\leq \sqrt{\frac{1}{2}}$, together with Lemma \ref{mini lemma}(i), to
show that $\mathbf{T}$ is $k$-hyponormal. \ To calculate $\sigma _{T}(%
\mathbf{T})$ and $\sigma _{Te}(\mathbf{T})$, we use Lemma \ref{lem1}, the
projection property for the Taylor spectrum, and direct verification of the
exactness of the Koszul complex associated with $\mathbf{T}$ at the middle
stage, using the technique in \cite[Corollary 4.3(ii)]{CuSa}. \qed

\begin{remark}
Careful analysis of the proof of Theorem \ref{khypo of ex} reveals that for $%
k\geq 1$ and $y_{0}\leq \sqrt{\frac{1}{2}}$, $\mathbf{T}$ is $k$-hyponormal.
\ Thus, $\mathbf{T}$\ is indeed subnormal in the interval $(0,\sqrt{\frac{1}{%
2}}]$, and its Berger measure is $\mu =y_0^2 \delta _{1}\times \delta
_{1}+(\xi _{\kappa}-y_0^2 \delta _{1})\times \delta _{0}$.
\end{remark}




\end{document}